\documentstyle{amsppt}
\magnification=\magstep1
\pagewidth{5.2in}
\pageheight{7.6in}
\abovedisplayskip=6pt
\belowdisplayskip=6pt

\topmatter
\title On the volume conjecture for hyperbolic knots \endtitle
\author Yoshiyuki Yokota
\endauthor
\affil
Graduate School of Mathematics, Kyushu University\\
Fukuoka, 812-8581, Japan\\
e-mail: jojo\@math.kyushu-u.ac.jp\\
From October, 2000: Department of Mathematics, Tokyo Metropolitan University\\
Tokyo, 192-0397, Japan
\endaffil
\email jojo\@math.kyushu-u.ac.jp \endemail
\endtopmatter

\head 1. Introduction \endhead

In [1],
R.M. Kashaev introduced certain invariants of oriented links
motivated by his study of quantum dilogarithm functions.
Since the classical dilogarithm functions
are related to the hyperbolic volumes,
he naturally expected that, for hyperbolic knots,
the asymptotic behaviors of his invariants determine their volumes,
which is in fact confirmed for a few hyperbolic knots by himself in [2].
Some other examples are now given in [5],
where the asymptotic behaviors of the invariants
suggestively determine not only the volumes but also the Chern-Simons invariants.

Later, in [4],
H. Murakami and J. Murakami have shown that
Kashaev's invariant coincides with
certain special value of the colored Jones function,
and then reformulated Kashaev's conjecture, that is,
the asymptotic behavior of the colored Jones function determines
the simplicial volume of a knot,
which is now called the {\it volume conjecture} of knots.
For torus knots, known as typical non-hyperbolic knots,
this conjecture is recently confirmed
by Kashaev and O. Tirkkonen in [3].

The purpose of this article is
to give a rough, and so not yet complete unfortunately, proof of Kashaev's conjecture,
that is, the volume conjecture for hyperbolic knots.
In fact, for a hyperbolic knot in $S^3$,
the quantum factorials in its Kashaev's invariant
naturally correspond to
the tetrahedra in an ideal triangulation of its complement,
which establish a surprising correspondence
between the stationary phase equations for the invariant
and the hyperbolicity equations for the triangulation.
Then, the asymptotic behavior of the invariant is determined by the promised solution,
which is nothing but the volume of the knot complement
because the quatum factorials asymptotically goes to the dilogarithm functions.

The author is grateful to Hirotaka Akiyoshi, Hitoshi Murakami
and Jun Murakami for stimulating discussions.
He is also grateful to all the participants
in the meeting \lq\lq Volume conjecture" in October, 1999
at International Institute for Advanced Studies, Kyoto.

\head 2. Ideal triangulations \endhead

Throughout of this article, $M$ denotes the complement of a hyperbolic knot $K$ in $S^3$.
In this section, associated with a diagram of $K$,
we construct an ideal triangulation of $M$
from which the hyperbolicity equations for $M$ follows quite nicely.
We suppose $K$ does not hit the two poles $\pm\infty$ of $S^3$
and denote $M$ with $\pm\infty$ removed by $\dot M$.

Let $D$ be a diagram of $K$ in $S^2$ with $n$ crossings
and $R_0,\dots,R_{n+1}$ its faces.
Then, $D$ together with its dual graph gives a decomposition $\Cal D$ of $S^2$
into $n$ octagons $Q_1,\dots,Q_n$,
where each $Q_\nu$ is further divided into four quadrangles as shown in Figure 1.
We put $\Cal Q_\nu=\{\mu|Q_\nu\cap R_\mu\ne\emptyset\}$ and
$\Cal R_\mu=\{\nu|Q_\nu\cap R_\mu\ne\emptyset\}$,
where we can suppose $|\Cal Q_\nu|=4$ for any $\nu$,
otherwise we can reduce the number of crossings of $D$.

$$\boxed{\text{Figure 1}}$$

Due to D. Thurston [6],
$\dot M$ decomposes into ideal octahedra $P_1,\dots,P_n$,
where each $P_\nu$ corresponds to $Q_\nu$
and further decomposes into four ideal tetrahedra as shown in Figure 2.
Thus, we have an ideal traiangulation $\dot\Cal S$ of $\dot M$
which associates $Q_\nu\cap R_\mu$ of $\Cal D$
with an ideal tetrahedron $S_{\nu\mu}$.
Notice that four tetrahedra corresponding to $\Cal Q_\nu$ share an edge $\dot E_\nu$
while all the tetrahedra corresponding to $\Cal R_\mu$ also share an edge $\dot F_\mu$.
As usual, we put a hyperbolic structure on $S_{\nu\mu}$
by assigning a complex number $z_{\nu\mu}$
to a pair of opposite edges of $S_{\nu\mu}$
corresponding to $\dot E_\nu$ and $\dot F_\mu$.

$$\boxed{\text{Figure 2}}$$

To make $\dot\Cal S$ an ideal triangulation of $M$,
we have to specify a non-singular point on $D$,
where we meet an overpass $A$ followed by an undercrossing $X$ in one side
and an underpass $B$ followed by an overcrossing $Y$ in the other side.
Without loss of generality, we can suppose $\Cal Q_1\cap\Cal Q_n=\{0,n+1\}$
and that $A,X,Y,B$ are covered by
$\cup_{\nu=1}^aQ_\nu,Q_x,Q_y,\cup_{\nu=b}^nQ_\nu$ respectively,
where $a<x\le y<b$,
otherwise we can reduce $n$, the number of crossings in $D$.
For simplicity, we put
$\Cal B=\{1,\dots,a,b,\dots,n\}$,
$\Cal R=\Cal R_0\cup\Cal R_{n+1}$
and $\Cal Q=\cup_{\nu\in\Cal B}\Cal Q_\nu$.
Then, we can further suppose $\Cal B\cap\Cal R=\{1,n\}$
and $\Cal Q_a\cap\Cal Q_x\cap\Cal Q_y\cap\Cal Q_b=\emptyset$,
and so $x\ne y$ in particular,
otherwise we can reduce $n$ again
or change the specified point to $D\cap Q_a\cap Q_x$ and so on,
which should stop sometime because $K$ can not be a satellite of an elementary torus link.

$$\boxed{\text{Figure 3}}$$

Now, we contract the ideal bigon $P_1\cap P_n$
bounded by $L=\{\pm\infty\}\cup\dot F_0\cup\dot F_{n+1}$,
which makes the ideal tetrahedra touching the bigon degenerate,
and come up with an ideal triangulation $\Cal S$ of $M$.
We choose a cusp cross section $T$ of $M$ so that $T$ does not touch the bigon,
and denote by $\dot\Cal T$ and $\Cal T$
the triangulations of $T$ induced by $\dot\Cal S$ and $\Cal S$ respectively.
Then, $(\dot\Cal S,\dot\Cal T)$ and $(\Cal S,\Cal T)$ are related as follows.

If $\nu\in\{1,n\}$, then $P_\nu$ intersects the bigon in two triangles
and so is truncated by $T$ as shown in Figure 4a.
Thus, no tetrahedron in $P_\nu$ survives in $\Cal S$ as shown in Figure 4b.
In particular, $S_{\nu0}$ and $S_{\nu(n+1)}$ degenerate into an edge in $\Cal S$.
Note that $\Cal Q_\nu\setminus\Cal Q=\emptyset$.

$$\boxed{\text{Figure 4}}$$

If $\nu\in\Cal B\setminus\{1,n\}$,
then two edges of $P_\nu$ are identified with an edge in the bigon,
and so is truncated by $T$ as shown in Figure 5a.
Thus, no tetrahedron in $P_\nu$ survives in $\Cal S$ as shown in Figure 5b,
and $\Cal Q_\nu\setminus\Cal Q=\emptyset$ again.

$$\boxed{\text{Figure 5}}$$

If $\nu\in\{x,y\}\cap\Cal R$,
then $P_\nu$ intersects the bigon in two edges,
and so is truncated by $T$ as shown in Figure 6a.
Thus, one tetrahedron, $S_{\nu c}$ say, in $P_\nu$ survives in $\Cal S$
as shown in Figure 6b.
Note that $\Cal Q_\nu\setminus\Cal Q=\{c\}$.

$$\boxed{\text{Figure 6}}$$

If $\nu\in\{x,y\}\setminus\Cal R$,
then $P_\nu$ intersects the bigon in one edge,
and so is truncated by $T$ as shown in Figure 7a.
Thus, two tetrahedra, $S_{\nu e},S_{\nu f}$ say, in $P_\nu$ survive in $\Cal S$.
Note that $\Cal Q_\nu\setminus\Cal Q=\{e,f\}$.

$$\boxed{\text{Figure 7}}$$

If $\nu\in\Cal R\setminus\{1,x,y,n\}$,
then $P_\nu$ intersects the bigon in one edge,
and so is truncated by $T$ as shown in Figure 8a.
Thus, three tetrahedra other than $S_{\nu g}$, where $g$ is 0 or $n+1$,
in $P_\nu$ survive in $\Cal S$ as shown in Figure 8b.
Note that $\Cal Q_\nu\setminus\Cal Q=\Cal Q_\nu\setminus\{g\}$.

$$\boxed{\text{Figure 8}}$$

Otherwise, that is, if $\nu\not\in\Cal B\cup\{x,y\}\cup\Cal R$,
$P_\nu$ does not touch the bigon, and so no degeneration occurs as shown in Figure 9.
Note that $\Cal Q_\nu\setminus\Cal Q=\Cal Q_\nu$.

$$\boxed{\text{Figure 9}}$$

Such observations enable us to write down $\Cal T$
and so the hyperbolicity equations for $M$ explicitly.
In what follows, we illustrate how to do this
through an example depicted in Figure 10.

$$\boxed{\text{Figure 10}}$$

Let $N(K),N(L),N(\pm\infty)$ be the regular neighborhoods
of $K,L,\pm\infty$ in $S^3$,
and $\dot\Cal K,\dot\Cal L,\Cal D_\pm$ the decompositions
of $\partial N(K),\partial N(L),\partial N(\pm\infty)$
induced by $\dot\Cal S$ respectively.
Then, it is not difficult to write down $\dot\Cal K$
which associates an edge of $D$ with either an annulus or a pinched annulus,
divided into four triangles anyway,
according as the edge is alternating or not.

$$\boxed{\text{Figure 11}}$$

On the other hand, $\Cal D_+$ is $\Cal D$ with the overpasses of $D$ contracted,
while $\Cal D_-$ is the mirror image of $\Cal D$ with the underpasses of $D$ contracted,
and $\dot\Cal L$ is then given by gluing $\Cal D_+$ and $\Cal D_-$
along $\dot F_0$ and $\dot F_{n+1}$.

$$\boxed{\text{Figure 12}}$$

In $S^3$ with $N(K)\cup N(L)$ removed,
the bigon becomes a properly embedded annulus 
and so $\dot\Cal T$ is obtained by gluing $\dot\Cal K$ and $\dot\Cal L$ along the annulus
as shown in Figure 13.
Then, by contracting a number of edges,
the dotted ones in Figure 13, in $\dot\Cal T$ as we have observed above,
we obtain $\Cal T$
which naturally decomposes into two triangulations $\Cal K$ and $\Cal L$ of annuli
corresponding to $\dot\Cal K$ and $\dot\Cal L$.
Note that $\Cal K\cap\Cal L$ represents two meridians of $K$
each of which consists of two edges.

$$\boxed{\text{Figure 13}}$$

To describe the hyperbolicity equations for $M$,
we consider a diagram $G$ with two trivalent vertices
which is obtained from $D$ by removing the simple arc between $X$ and $Y$.
Let $m$ denote the number of crossings of $G$,
and $\Cal M_0\cup\dots\cup\Cal M_{m+1}$ a partition of $\{1,\dots,n\}$
so that $\cup_{\mu\in\Cal M_\lambda}R_\mu$ form a face of $G$
other than $R_0\cup R_{n+1}$,
where we can suppose
$\Cal Q_a\cap\Cal Q_x\subset\Cal M_0$,
$\Cal Q_b\cap\Cal Q_y\subset\Cal M_{m+1}$.
Let $\Cal E$ denote the sets of edges of $G$
and $\Cal F$ the set of edges of $G$ lying in $\partial R_0\cup\partial R_{n+1}$.

$$\boxed{\text{Figure 14}}$$

If $\nu\not\in\Cal B\cup\{x,y\}$,
$\dot E_\nu$ survives alone as an edge $E_\nu$ in $\Cal S$,
which connects two vertices in $\Cal K\setminus\Cal L$,
and the edge relation around $E_\nu$ reads
$$\prod_{\mu\in\Cal Q_\nu\setminus\Cal Q}z_{\nu\mu}=1.$$
On the other hand, if $\lambda\not\in\{0,m+1\}$,
the edges in $\{\dot F_\mu|\mu\in\Cal M_\lambda\}$
reduce to an edge $F_\lambda$ in $\Cal S$,
which connects two vertices in $\Cal L\setminus\Cal K$,
and the edge relation around $F_\lambda$ reads
$$\prod_{\mu\in\Cal M_\lambda}
\prod_{\nu\in\Cal R_\mu\setminus\Cal B}z_{\nu\mu}=1.$$
If $\nu\in\{x,y\}$ and $\lambda\in\{0,m+1\}$,
$\dot E_\nu$, together with some other edges in $\dot\Cal S$,
survives as an edge $E_\nu$ in $\Cal S$,
connecting two vertices in $\Cal K\setminus\Cal L$ and $\Cal K\cap\Cal L$,
while the edges in $\{\dot F_\mu|\mu\in\Cal M_\lambda\}$,
together with some other edges in $\dot\Cal S$ too,
reduce to an edge $F_\lambda$ in $\Cal S$,
connecting two vertices in $\Cal L\setminus\Cal K$ and $\Cal K\cap\Cal L$.
Thus, the edge relations around $E_x,E_y$ should be read from $\Cal K\setminus\Cal L$,
and those around $F_0,F_{m+1}$ may be substituted
with the cusp conditions along the two annuli lying on the borders of $\Cal K$ and $\Cal L$,
which read
$$\prod_{\mu\in\Cal M_0}\prod_{\nu\in\Cal R_\mu\setminus\Cal B}z_{\nu\mu}
=\prod_{\mu\in\Cal Q_x\setminus\Cal Q}z_{x\mu},\quad
\prod_{\mu\in\Cal M_{m+1}}\prod_{\nu\in\Cal R_\mu\setminus\Cal B}z_{\nu\mu}
=\prod_{\mu\in\Cal Q_y\setminus\Cal Q}z_{y\mu}.$$

Note that the set of solutions to these equations coincides with
the set of functions from $\Cal E$ to $\Bbb C$ which take 1 on $\Cal F$.
In fact, for such $z$,  $z_{\nu\mu}$ is given by
$z(\varphi_{\nu\mu})/z(\psi_{\nu\mu})$,
where $\varphi_{\nu\mu},\psi_{\nu\mu}\in\Cal E$ touch $Q_\nu\cap R_\mu$
as shown in Figure 15.
We here put $\epsilon(\nu,\mu)=1$ or $-1$
according as $\psi_{\nu\mu}$ is over $\varphi_{\nu\mu}$ in $D$ or not.

$$\boxed{\text{Figure 15}}$$

Now, any other equation can be read from $\Cal K\setminus\Cal L$.
Let $\dot\Cal A_\varphi$ be an annulus, which may be pinched, in $\dot\Cal K$
corresponding to $\varphi\in\Cal E$
and $\Cal A_\varphi$ denote $\dot\Cal A_\varphi$ in $\Cal K$.
Furthermore, $\Cal U$ denotes the set of vertices in $\Cal K$ correponding to $E_\nu$'s
and $\Cal V$ denotes the set of vertices in $\Cal K$
which does not lie in $\partial\Cal A_\varphi$ for any $\varphi\in\Cal E$.
If $\varphi\in\Cal F$,
then $\Cal A_\varphi$ is a circle or a pinched annulus in $\Cal K$
such that $\Cal A_\varphi\cap\Cal V=\emptyset$,
and so no equation arises.
We then suppose $\varphi\in\Cal E\setminus\Cal F$ is surrounded
by $\alpha,\beta,\gamma,\delta\in\Cal E$ as shown in Figure 16,
and put
$$
\epsilon(\varphi)=\cases
1 & \text{if $\varphi$ is over $\alpha\cup\beta$ and $\gamma\cup\delta$ in $D$,}\\
-1 & \text{if $\varphi$ is under $\alpha\cup\beta$ and $\gamma\cup\delta$ in $D$,}\\
0 & \text{otherwise.}\\
\endcases
$$

$$\boxed{\text{Figure 16}}$$

If $\epsilon(\varphi)=1$,
then $\Cal A_\varphi$ is a pinched annulus 
with $\Cal A_\varphi\cap\Cal V\ne\emptyset$ as shown in Figure 17,
and so the edge relation around $\Cal A_\varphi\cap\Cal V$ reads
$${1-z(\alpha)/z(\varphi)\over1-z(\beta)/z(\varphi)}\cdot
{1-z(\delta)/z(\varphi)\over1-z(\gamma)/z(\varphi)}=1,$$
where $z(\varepsilon)/z(\varphi)$ is deleted
if $\varepsilon\in\{\alpha,\beta,\gamma,\delta\}$ is empty.

$$\boxed{\text{Figure 17}}$$

If $\epsilon(\varphi)=-1$,
then $\Cal A_\varphi$ is a pinched annulus 
with $\Cal A_\varphi\cap\Cal V\ne\emptyset$ too as shown in Figure 18,
and so the edge relation around $\Cal A_\varphi\cap\Cal V$ reads
$${1-z(\varphi)/z(\beta)\over1-z(\varphi)/z(\alpha)}\cdot
{1-z(\varphi)/z(\gamma)\over1-z(\varphi)/z(\delta)}=1,$$
where $z(\varphi)/z(\varepsilon)$ is deleted
if $\varepsilon\in\{\alpha,\beta,\gamma,\delta\}$ is empty.

$$\boxed{\text{Figure 18}}$$

If $\epsilon(\varphi)=0$, we can suppose
$\varphi$ is over $\alpha\cup\beta$ and under $\gamma\cup\delta$ in $D$.
Then, $\Cal A_\varphi$ is an annulus in $\Cal K$ as shown in Figure 19,
and the cusp condition along $\Cal A_\varphi$ reads
$${1-z(\alpha)/z(\varphi)\over1-z(\beta)/z(\varphi)}
={1-z(\varphi)/z(\gamma)\over1-z(\varphi)/z(\delta)},$$
where $z(\varepsilon)/z(\varphi)$ is deleted
if $\varepsilon\in\{\alpha,\beta\}$ is empty
and $z(\varphi)/z(\varepsilon)$ is deleted
if $\varepsilon\in\{\gamma,\delta\}$ is empty.
The edge relations around the vertices not in $\Cal U\cup\Cal V$
may be substituted with these cusp conditions,
and so we obtain a complete system of hyperbolicity equations for $M$
which has $2m+3$ unknowns and $2m+3$ equations 
corresponding to $\Cal E\setminus\Cal F$.

$$\boxed{\text{Figure 19}}$$

To be more precise,
we shall introduce a potential function for these equations
by using Euler's dilogarithm function
$$
\text{Li}_2(\omega)=-\int_0^\omega{\log(1-w)\over w}dw.
$$
For $z:\Cal E\to\Bbb C$, we put
$$
V(z)=\sum_{\nu\not\in\Cal B}
\sum_{\mu\in\Cal Q_\nu\setminus\Cal Q}
V_{\nu\mu}(z)-2\pi\sqrt{-1}\sum_{\varphi\in\Cal E\setminus\Cal F}
\epsilon(\varphi)\log z(\varphi),
$$
where
$$
V_{\nu\mu}(z)=\epsilon(\nu,\mu)\cdot
\{\text{Li}_2(
z(\varphi_{\nu\mu})^{\epsilon(\nu,\mu)}/z(\psi_{\nu\mu})^{\epsilon(\nu,\mu)}
)-\pi^2/6\}.
$$
Then, a simple calculation shows
$$\text{Im}\,V(z)=\sum_{\nu\not\in\Cal B}
\sum_{\mu\in\Cal Q_\nu\setminus\Cal Q}
D(z(\varphi_{\nu\mu})/z(\psi_{\nu\mu}))
+\sum_{\varphi\in\Cal E\setminus\Cal F}
\log|z(\varphi)|\cdot\text{Im}\,z(\varphi){\partial V(z)\over\partial z(\varphi)},$$
where
$$D(\omega)=\text{Im}\,\text{Li}_2(\omega)+\log|\omega|\arg(1-\omega).$$
It is well-known that the hyperbolic volume of $S_{\nu\mu}$ is given by
$D(z(\varphi_{\nu\mu})/z(\psi_{\nu\mu}))$,
and so our observations can be summarized as follows.

\proclaim{Proposition 1}
The hyperbolicity equations for $M$ associated to $D$ are given by
$$\left\{\partial V_0(z)/\partial z(\varphi)=0\,|\,
\varphi\in\Cal E\setminus\Cal F\right\},$$
where $V_0(z)$ is a branch of $V(z)$,
and the volume $\text{\rm vol}(M)$ of $M$ is then given by
$$\text{\rm vol}(M)=\text{\rm Im}\,V_0(z_0),$$
where $z_0$ is the promised solution to the hyperbolicity equations.
\endproclaim

\head 3. Kashaev's invariants \endhead

Let $N$ be a positive integer, which will be sent to $\infty$,
and $\Cal N=\{0,1,\dots,N-1\}$.
For $h\in\Bbb Z$, we denote by $[h]\in\Cal N$ the residue modulo $N$.
In this section,
we compute Kashaev's invariant $\langle K\rangle_N$ of $K$,
which is nothing but the $N$-colored Jones function of $K$
evaluated at $q=\exp 2\pi\sqrt{-1}/N$ due to [4],
and detect its asymptotic behavior.
Put
$$
R^{ij}_{kl}=
{Nq^{-1-(k-j)(i-l+1)}\theta_{ijkl}
\over(\bar q)_{[i-j]}(q)_{[j-l]}(\bar q)_{[l-k-1]}(q)_{[k-i]}},
\quad
\bar R^{ij}_{kl}=
{Nq^{1+(i-l)(k-j+1)}\theta_{ijkl}
\over(q)_{[i-j]}(\bar q)_{[j-l]}(q)_{[l-k-1]}(\bar q)_{[k-i]}}
$$
for $i,j,k,l\in\Cal N$, where $(\omega)_{[h]}=(1-\omega)(1-\omega^2)\dots(1-\omega^{[h]})$
and
$$\theta_{ijkl}=\cases
1 & \text{if $[i-j]+[j-l]+[l-k-1]+[k-i]=N-1$,}\\
0 & \text{otherwise.}
\endcases$$
Furthermore, we write $k\in[i,j]$ if $[i-k]+[k-i]=[i-j]$.

In what follows,
we suppose $D$ is a closed braid,
where $A,X,Y,B$ are arranged in this order along its orientation,
and put $D$ in $\Bbb R^2$ by removing the specified point.
Then, we can suppose
$|\Cal Q_\nu|=4$, $\Cal B\cap\{x,y\}=\emptyset$ and $\Cal B\cap\Cal R=\{1,n\}$ as before,
otherwise we can reduce the number of strings of $D$.
Of course, we can suppose $\Cal Q_a\cap\Cal Q_x\cap\Cal Q_y\cap\Cal Q_b=\emptyset$.

Let $\Cal X$  be the set of maxima of $D$ where $D$ is oriented clockwise
and $\Cal Y$ the set of maxima  of $D$ where $D$ is oriented counter-clockwise.
In what follows,
such maxima and minima are regarded as vertices of $D$ together with crossings,
and the set of the edges of $D$ is then denoted by $\dot\Cal E$.
Furthermore, as shown in Figure 20,
$\alpha_\nu,\beta_\nu,\gamma_\nu,\delta_\nu$ denote
the four edges in $\dot\Cal E$ incident to the crossing with sign $\epsilon(\nu)$ in $Q_\nu$
while $\alpha_\xi,\delta_\xi$ and $\alpha_\eta,\delta_\eta$
denote the edges in $\dot\Cal E$ incident to $\xi\in\Cal X$ and $\eta\in\Cal Y$ respectively.

$$\boxed{\text{Figure 20}}$$

A state of $D$ is a function $\sigma:\dot\Cal E\to\Cal N$
which assigns 0 to the non-compact edges of $D$.
For such $\sigma$, we put
$$\langle D|\sigma\rangle_\nu=\cases
R_{\sigma(\gamma_\nu)\sigma(\delta_\nu)}^{\sigma(\alpha_\nu)\sigma(\beta_\nu)} &
\text{if $\epsilon(\nu)=+1$,}\\
\bar R_{\sigma(\delta_\nu)\sigma(\gamma_\nu)}^{\sigma(\beta_\nu)\sigma(\alpha_\nu)} &
\text{if $\epsilon(\nu)=-1$.}
\endcases$$
for $\nu\in\{1,\dots,n\}$.
Then, Kashaev's invariant $\langle K\rangle_N$ of $K$ is given by
$$\langle K\rangle_N=\prod_{\nu=1}^nq^{\epsilon(\nu)/2}
\sum_{\sigma\in\Cal Z}\langle D|\sigma\rangle,$$
where $\Cal Z$ denotes the set of states of $D$ and
$$\langle D|\sigma\rangle=\prod_{\nu=1}^n\langle D|\sigma\rangle_\nu
\prod_{\xi\in\Cal X}
\{-q^{1/2}\delta_{\sigma(\alpha_\xi)+1,\sigma(\delta_\xi)}\}
\prod_{\eta\in\Cal Y}
\{-q^{-{1/2}}\delta_{\sigma(\alpha_\eta)+1,\sigma(\delta_\eta)}\}.$$
However, a lot of states do not contribute to $\langle K\rangle_N$ at all.
To describe them, in what follows, we denote
$$\{[\epsilon(\nu)\sigma(\alpha_\nu-\beta_\nu)],
[\epsilon(\nu)\sigma(\beta_\nu-\delta_\nu)],
[\epsilon(\nu)\sigma(\gamma_\nu-\alpha_\nu)],
[\epsilon(\nu)\sigma(\delta_\nu-\gamma_\nu)-1]\}$$
by $\{\sigma(\nu,\mu)\,|\,\mu\in\Cal Q_\nu\}$ for $\nu\in\{1,\dots,n\}$,
where $\sigma(\varphi-\psi)$ stands for $\sigma(\varphi)-\sigma(\psi)$.

\proclaim{Lemma 2} If $\langle D|\sigma\rangle\ne0$,
$$\sum_{\mu\in\Cal Q_\nu}\sigma(\nu,\mu)
=\sum_{\nu\in\Cal R_\mu}\sigma(\nu,\mu)=N-1.$$
Furthermore, $\sigma(\nu,0)=0$ for $\nu\in\Cal R_0$
and $\sigma(\nu,n+1)=0$ for $\nu\in\Cal R_{n+1}$.
\endproclaim

\demo{Proof} By definition, we have
$$\sum_{\mu\in\Cal Q_\nu}\sigma(\nu,\mu)=N-1,\quad
\sum_{\nu=1}^n\sum_{\mu\in\Cal Q_\nu}\sigma(\nu,\mu)
=\sum_{\mu=0}^{n+1}\sum_{\nu\in\Cal R_\mu}\sigma(\nu,\mu)$$
and so
$$\sum_{\mu=0}^{n+1}\sum_{\nu\in\Cal R_\mu}\sigma(\nu,\mu)=n(N-1).$$
On the other hand, we can observe
$$\sum_{\nu\in\Cal R_\mu}\sigma(\nu,\mu)\ge N-1$$
unless $\mu\in\{0,n+1\}$, and so Lemma 2 follows immediately. \qed
\enddemo

Lemma 2 definitely reduces the definition of $\langle K\rangle_N$ but not enough.
We can further reduce certain $q$-factorials in $\langle K\rangle_N$
by using

\proclaim{Lemma 3} We have
$$\gather
\sum_{i\in[k,j]}q^{-i}\bar R_{kl}^{ij}=\delta_{j,k}q^{1-l},\ 
\sum_{j\in[i,l]}q^{-j}R_{kl}^{ij}=\delta_{i,l}q^{-1-k},\\
\sum_{k\in[l,i]}q^k\bar R_{kl}^{ij}=\delta_{i+1,l}q^j,\ 
\sum_{l\in[j,k]}q^lR_{kl}^{ij}=\delta_{j,k+1}q^i\ 
\endgather$$
and
$$\gather
\sum_{i\in[k,j]}q^{-i}R_{kl}^{ij}
={-Nq^{-1-k}\over(\bar q)_{[j-l]}(q)_{[l-k-1]}},\ 
\sum_{j\in[i,l]}q^{-j}\bar R_{kl}^{ij}
={-Nq^{1-l}\over(\bar q)_{[l-k-1]}(q)_{[k-i]}},\\
\sum_{k\in[j,l]}q^kR_{kl}^{ij}
={-Nq^{-1+i}\over(q)_{[i-j]}(\bar q)_{[j-l]}},\ 
\sum_{l\in[k,i]}q^l\bar R_{kl}^{ij}
={-Nq^{1+j}\over(\bar q)_{[i-j]}(q)_{[k-i]}}.
\endgather$$
\endproclaim

\demo{Proof} We shall prove
$$\sum_{i\in[k,j]}q^{-i}\bar R_{kl}^{ij}=\delta_{j,k}q^{1-l},\ 
\sum_{i\in[k,j]}q^{-i}R_{kl}^{ij}
={-Nq^{-1-k}\over(\bar q)_{[j-l]}(q)_{[l-k-1]}}.$$
By Lemma A.1 of [4], we have
$$\sum_{k\in[i,j]}
{q^{-k(i-j)}\over(q)_{[i-k]}(\bar q)_{[k-j]}}=\delta_{i,j},\quad
\sum_{k\in[i,j]}
{q^{-k(i-j+1)}\over(q)_{[i-k]}(\bar q)_{[k-j]}}=(-1)^{i-j}q^{-(i+j)(i-j+1)/2},$$
and so
$$
\sum_{i\in[k,j]}q^{-i}\bar R_{kl}^{ij}
=\sum_{i\in[k,j]}
{q^{-i(j-k)}\over(q)_{[i-j]}(\bar q)_{[k-i]}}
{Nq^{1+l(j-k-1)}\over(q)_{[l-k-1]}(\bar q)_{[j-l]}}=\delta_{j,k}q^{1-l},
$$

\newpage

$$
\sum_{i\in[k,j]}q^{-i}R_{kl}^{ij}
=\sum_{i\in[k,j]}
{q^{-i(k-j+1)}\over(q)_{[k-i]}(\bar q)_{[i-j]}}
{Nq^{-1-(j-k)(l-1)}\over(q)_{[j-l]}(\bar q)_{[l-k-1]}}
={-Nq^{-1-k}\over(\bar q)_{[j-l]}(q)_{[l-k-1]}}.\ \qed
$$
\enddemo

Let $[\sigma]$ denote the set of states of $D$
which coincide with $\sigma$ on the edges lying in $G$.
Then, by Lemma 3,
$$
\align
\sum_{\tau\in[\sigma]}
\prod_{\nu\in\Cal B\cup\{x,y\}}&\langle D|\tau\rangle_\nu
=\prod_{\nu=1}^a\delta_{\sigma(\alpha_\nu),\sigma(\delta_\nu)}
\prod_{\nu=b}^n\delta_{\sigma(\beta_\nu),\sigma(\gamma_\nu)+\epsilon(\nu)}
\prod_{\nu\in\Cal B}q^{-\epsilon(\nu)}\\
\times&
{N^2\over
(q^{\epsilon(x)})_{[\epsilon(x)\sigma(\delta_x-\gamma_x)-1]}
(q^{\epsilon(y)})_{[\epsilon(y)\sigma(\alpha_y-\beta_y)]}}
\prod_{\nu\in\{x,y\}}
{q^{-\epsilon(\nu)}\over
(q^{-\epsilon(\nu)})_{[\epsilon_\nu\sigma(\beta_\nu-\delta_\nu)]}},
\endalign
$$
and so we come up with the following definition.
A state $\sigma$ of $D$ is said to be simple
if $\sigma(\varphi)=0$ for $\varphi\subset D\setminus G$ and
$$
\prod_{\nu=1}^a\delta_{\sigma(\alpha_\nu),\sigma(\delta_\nu)}
\prod_{\nu=b}^n\delta_{\sigma(\beta_\nu),\sigma(\gamma_\nu)+\epsilon(\nu)}
\prod_{\xi\in\Cal X}\delta_{\sigma(\alpha_\xi)+1,\sigma(\delta_\xi)}
\prod_{\eta\in\Cal Y}\delta_{\sigma(\alpha_\eta)+1,\sigma(\delta_\eta)}
\ne0.
$$
Thus, a simple state can be considered as a function from $\Cal E$ to $\Cal N$ and
$$
\align
\langle D|\sigma\rangle=&(-q^{1/2})^{|\Cal X|-|\Cal Y|}
\prod_{\nu\in\Cal B}q^{-\epsilon(\nu)}\times
\prod_{\nu\not\in\Cal B\cup\{x,y\}}\langle D|\sigma\rangle_\nu\\
\times&
{N^2\over
(q^{\epsilon(x)})_{[\epsilon(x)\sigma(\delta_x-\gamma_x)-1]}
(q^{\epsilon(y)})_{[\epsilon(y)\sigma(\alpha_y-\beta_y)]}}
\prod_{\nu\in\{x,y\}}
{q^{-\epsilon(\nu)}\over
(q^{-\epsilon(\nu)})_{[\epsilon(\nu)\sigma(\beta_\nu-\delta_\nu)]}}.
\endalign
$$
for a simple state $\sigma$.
Furthermore, by using
$$
(q)_{[h]}=\pm(-1)^hq^{h(h+1)/2}(\bar q)_{[h]},
$$
we have
$$
R_{ij}^{kl}=
{\pm Nq^{i-k-1}\theta_{ijkl}\over(q)_{[i-j]}(\bar q)_{[j-l]}(q)_{[l-k-1]}(\bar q)_{[k-i]}},
\quad
\bar R_{ij}^{kl}=
{\pm Nq^{j-l+1}\theta_{ijkl}\over(\bar q)_{[i-j]}(q)_{[j-l]}(\bar q)_{[l-k-1]}(q)_{[k-i]}}
$$
and
$$
\langle D|\sigma\rangle_\nu
={-Nq^{\sigma(\alpha_\nu-\gamma_\nu)-\epsilon(\nu)}\over
\prod_{\mu\in\Cal Q_\nu\setminus\Cal Q}(q^{\epsilon(\nu,\mu)})_{\sigma(\nu,\mu)}}.
$$
Consequently, $\langle K\rangle_N$ can be rewritten as
$$
(-N)^{n-a-b}(-q^{1/2})^{|\Cal X|-|\Cal Y|}
\prod_{\nu=1}^nq^{-\epsilon(\nu)/2}
\sum_{\sigma\in\Cal Z_0}
\prod_{\nu\not\in\Cal B}
{q^{-(N-1)\,\sigma(\alpha_\nu-\gamma_\nu)}\over
\prod_{\mu\in\Cal Q_\nu\setminus\Cal Q}(q^{\epsilon(\nu,\mu)})_{\sigma(\nu,\mu)}},
$$
where $\Cal Z_0$ denotes the set of simple states of $D$.
It should be noted that the $q$-factorials above curiously corresponds to the tetrahedra in $\Cal S$.

%\head 4. Analysis \endhead

From now onward, we suppose $N$ is sufficiently large so that
$$
{1\over(q^{\pm1})_{[h]}}\sim
\exp{N\{\pm\text{Li}_2(q^{\pm[h]})\mp\pi^2/6\}\over2\pi\sqrt{-1}}.
$$
Then, due to Kashaev [2],
$\langle K\rangle_N$ can be approximated by the integral
$$
\int\dots\int \exp{NV(z)\over2\pi\sqrt{-1}}
\prod_{\varphi\in\Cal E\setminus\Cal F}dz(\varphi),
$$
where $\{z(\varphi)\,|\,\varphi\in\Cal E\setminus\Cal F\}$ correspond to
$\{q^{\sigma(\varphi)}\,|\,\varphi\in\Cal E\setminus\Cal F\}$ in $\langle K\rangle_N$.
The maximal contributions to this integral comes from
the solutions to the stationary phase equations for the branches of $V(z)$
which of course contain the solution $z_0$ to
$$
\{\partial V_0(z)/\partial z(\varphi)=0\,|\,\varphi\in\Cal E\setminus\Cal F\},
$$
the hyperbolicity equations before.
Let $z_1$ be such a solution other than $z_0$
which is derived from the stationary phase equations for a branch $V_1(z)$ of $V(z)$.
Then, the following assumption is likely to be true,
see p.\,367 in [7],
but the author has no proof.

\demo{Assumption} $\text{Im}\,V_1(z_1)<\text{Im}\,V_0(z_0)$. \enddemo

If this is true,
the integral is simply approximated by
$$
\exp{NV_0(z_0)\over2\pi\sqrt{-1}}
$$
by the saddle point method, where we have to use another assumption.

\demo{Assumption} There is a deformation of the contour
to apply the saddle point method. \enddemo

Then, the absolute value of the invariant exponentially grows like
$$
e^{{N\over2\pi}\text{Im}\,V_0(z_0)}=e^{{N\over2\pi}\text{vol}(M)},
$$
by Proposition 1,
which proves the volume conjecture for hyperbolic knots.

\end{document}